\documentstyle[12pt]{article}
\def\Bbb R{{\rm \bf R}}
\def\proclaim#1{\vskip2mm{\bf #1}\em}
\def\endproclaim{\em \vskip2mm}
\def\tag#1{\eqno(#1)}
\def\gathered{\begin{array}{c}}
\def\endgathered{\end{array}}
\def\text{\mbox}

\begin{document}

\title {The enclosure method for inverse obstacle scattering over a finite time interval: V.
Using time-reversal invariance}
\author{Masaru IKEHATA\footnote{
Laboratory of Mathematics,
Graduate School of Engineering,
Hiroshima University,
Higashihiroshima 739-8527, JAPAN}}
\maketitle

\begin{abstract}
The wave equation is time-reversal invariant. The enclosure method using
a Neumann data generated by this invariance is introduced.  The method yields
the minimum ball that is centered at a given arbitrary point and encloses an unknown obstacle embedded in a known bounded domain
from a {\it single point} on the graph of the so-called response operator on
the boundary of the domain over a finite time interval.  The occurrence of the lacuna
in the solution of the free space wave equation is positively used.

\noindent
AMS: 35R30, 35L05

\noindent KEY WORDS: enclosure method, time-reversal invariance, inverse obstacle problem, wave equation,
non-destructive testing.
\end{abstract}


\section{Introduction}

The so-called inverse obstacle problem is a typical problem in the inverse problems community
and the solution has several possibilities of applications to non-destructive testing, sonar, radar, to name a few. 
See \cite{Is} for a survey about the uniqueness and stability issue.

This paper is concerned with the reconstruction or extraction issue, in particular, its methodology.
Succeeding to the previous studies about the {\it time domain enclosure method} for inverse obstacle problems governed by
the wave equation developed in \cite{IW00, IEO2, IEO3}, we further continue to pursue various possibilities of the method itself.
In \cite{EIV} the author has introduced a new version of the time domain enclosure method for inverse obstacle scattering problems
using the wave equation in a bounded domain over a finite time interval.
The method employs the Neumann data generated by taking the normal derivative of a solution of the wave equation 
in the whole space, on the boundary of the domain and yields the maximum ball that is centered at an arbitrary given point {\it outside} the domain,
and its exterior encloses an unknown obstacle embedded in the domain. The point is: it makes use
of a {\it single point} on the graph of the response operator associated with the wave equation in the domain.
The aim of this paper is to add one more point to this new version of the enclosure method.
It is a combination of the {\it time-reversal invariance} of the wave equation and the new
version of the enclosure method.  As a result one can find the minimum ball, with a fixed
center, enclosing the obstacle by using only single boundary measurement over a finite
time interval.

To clarify the essence of the idea we consider the same problem as in \cite{EIV}.

First let us recall the problem.
Let $\Omega$ be a bounded domain of $\Bbb R^3$ with $C^2$-boundary.
Let $D$ be a nonempty bounded open set of $\Bbb R^3$
with $C^2$-boundary such that $\overline D\subset\Omega$ and $\Omega\setminus\overline D$ is connected.

Given an arbitrary positive number $T$ and $f=f(x,t)$, $(x,t)\in\partial\Omega\times\,[0,\,T]$,
let $u=u_{f}(x,t)$, with $(x,t)\in(\Omega\setminus\overline D)\times\,[0,\,T]$, denote the solution of the following initial
boundary value problem for the classical wave equation:
$$\left\{
\begin{array}{ll}
\displaystyle
(\partial_t^2-\Delta)u=0 & \text{in $(\Omega\setminus\overline D)\times\,]0,\,T[$,}
\\
\\
\displaystyle
u(x,0)=0 & \text{in $\Omega\setminus\overline D$,}\\
\\
\displaystyle
\partial_tu(x,0)=0 & \text{in $\Omega\setminus\overline D$,}\\
\\
\displaystyle
\frac{\partial u}{\partial\nu}=0 & \text{on $\partial D\times\,]0,\,T[$,}\\
\\
\displaystyle
\frac{\partial u}{\partial\nu}=f(x,t) & \text{on $\partial\Omega\times\,]0,\,T[$.}
\end{array}
\right.
\tag {1.1}
$$
The problem considered in \cite{EIV} is the following.

{\bf\noindent Problem.}  Fix $T$ (to be determined later).  Assume that $D$ is unknown.
Find a suitable Neumann data $f$ in such a way that the wave
$u_f$ on $\partial\Omega$ over the time interval $[0,\,T]$ yields
information about the geometry of $D$.

What we found therein is: if $f$ is given by the normal derivative of a special solution of the Cauchy problem
for the classical wave equation
in $\Bbb R^3\times\,]0,\,T[$ with special initial data supported on an arbitrary fixed ball
outside $\Omega$, then one can extract the distance of the ball
to $D$ provided, roughly speaking,  $T$ is large enough.

In this paper, we give another choice of the Neumann data
that yields another information about the geometry of $D$.

Let $B$ be an open ball centered at $p\in\Bbb R^3$ with radius $\eta$
and denote by $\chi_B$ its characteristic function.
Define
$$\begin{array}{ll}
\displaystyle
\Psi_B(x)=(\eta-\vert x-p\vert)\chi_B(x), & x\in\Bbb R^3
\end{array}
$$
This function belongs to $H^1(\Bbb R^3)$ and $\text{supp}\,\Psi_B=\overline B$.
Unlike \cite{EIV}, in this paper we do not make a restriction on the position of $B$ relative to $\Omega$.

Let $v=v(x,t)$ be the solution of the following Cauchy problem for the classical wave equation:
$$
\left\{
\begin{array}{ll}
\displaystyle
(\partial_t^2-\Delta)v=0, & x\in\Bbb R^3, 0<t<T
\\
\\
\displaystyle
v(x,0)=0,  & x\in\Bbb R^3,\\
\\
\displaystyle
\partial_tv(x,0)=\Psi_B(x), & x\in\Bbb R^3.
\end{array}
\right.
\tag {1.2}
$$
It is well known that the solution $v$ takes the form
$$\displaystyle
v(x,t)=\frac{1}{4\pi t}\int_{\partial B_t(x)}\Psi_B(y)dS_y,
\tag {1.3}
$$
where
$$\displaystyle
B_t(x)=\{y\in\Bbb R^3\,\vert\,\vert y-x\vert<t\}.
$$
From the form of (1.3), we see that 
$$
\displaystyle
\text{supp}\,v(\,\cdot\,,T)\cup\text{supp}\,\partial_tv(\,\cdot\,,T)
\subset \overline{B_{T+\eta}(p)}
\tag {1.4}
$$
and
$$
\displaystyle
\text{supp}\,v(\,\cdot\,,T)\cup\text{supp}\,\partial_tv(\,\cdot\,,T)
\subset \Bbb R^3\setminus B_{T-\eta}(p),
\tag {1.5}
$$
where $B_{T\pm\eta}(p)=\{x\in\Bbb R^3\,\vert\,\vert x-p\vert<T\pm\eta\}$.
In \cite{EIV} we made use of (1.4) only, however, in this paper
we make use of also property (1.5), which is a quantitative expression of
{\it occurrence of lacuna} (cf. \cite{Du}).
It is a character of the wave equation in {\it odd} dimensions.

In this paper, we always choose $T$ in such a way that
$$\displaystyle
\Omega\subset B_{T-\eta}(p),
$$
that is,
$$\displaystyle
T-\eta\ge R_{\Omega}(p),
\tag {1.6}
$$
where
$$\displaystyle
R_{\Omega}(p)=\sup_{x\in\Omega}\vert x-p\vert.
$$
Define
$$\begin{array}{lll}
\displaystyle
f_{B,T}(x,t)=\frac{\partial}{\partial\nu}v(x,T-t), & x\in\partial\Omega, & 0\le t\le T.
\end{array}
\tag {1.7}
$$
This is the special $f$ mentioned above.
Note that property (1.5) and the {\it time-reversal invariance} of the wave equation yield
the function $v^*(x,t)=v(x,T-t)$ for $x\in\Omega$ and $0<t<T$ satisfies (1.1) with
$D=\emptyset$ and $f=f_{B,T}$.  Then a combination of a standard lifting argument 
and the theory of $C_0$-semigroups \cite{Y}
enables us to solve (1.1), with $f=f_{B,T}$ uniquely in the class
$$\displaystyle
C^2([0,\,T], L^2(\Omega\setminus\overline D))
\cap
C^1([0,\,T], H^1(\Omega\setminus\overline D))\cap C([0,\,T], H^2(\Omega\setminus\overline D)).
$$
See \cite{EIV} for this argument and \cite{I} for the solvability of the reduced problems,
which are initial boundary value problems for hyperbolic equations with {\it homogeneous} boundary conditions.

Having the solution $u=u_f$ of (1.1) with $f=f_{B,T}$ given by (1.7),
set
$$\begin{array}{lll}
\displaystyle
w_{B,T}(x)=w_{B,T}(x,\tau)=\int_0^T e^{-\tau t} u_f(x,t)dt, & x\in\Omega\setminus\overline D, & \tau>0,
\end{array}
\tag {1.8}
$$
and
$$\begin{array}{lll}
\displaystyle
w_{B,T}^*(x)=w_{B,T}^*(x,\tau)=\int_0^T e^{-\tau t} v(x,T-t)dt, & x\in\Bbb R^3, & \tau>0.
\end{array}
\tag {1.9}
$$

Define the {\it indicator function}
$$\begin{array}{ll}
\displaystyle
I_{\partial\Omega}(\tau;B,T)
=\int_{\partial\Omega}(w_{B,T}-w_{B,T}^*)\frac{\partial w_{B,T}^*}{\partial\nu}\,dS,
&
\tau >0.
\end{array}
$$
Define
$$\displaystyle
R_D(p)=\sup_{x\in D}\,\vert x-p\vert.
$$
Note that again, in this paper, we always assume that $T$ satisfies (1.6).

\proclaim{\noindent Theorem 1.1.}
(i)  Let $\eta$ satisfy 
$$\displaystyle
\eta+2R_D(p)>R_{\Omega}(p).
\tag {1.10}
$$
Then, there exists a positive number $\tau_0$ such that $I_{\partial\Omega}(\tau;B,T)>0$
for all $\tau\ge\tau_0$, and
we have
$$\displaystyle
\lim_{\tau\rightarrow\infty}
\frac{1}{\tau}\log I_{\partial\Omega}(\tau;B,T)
=-2\left\{(T-\eta)-R_D(p)\right\}.
\tag {1.11}
$$

(ii) If  $T>2\{(T-\eta)-R_D(p)\}$, then
$$\displaystyle
\lim_{\tau\rightarrow\infty}
e^{\tau T}I_{\partial\Omega}(\tau;B,T)
=\infty.
$$

(iii)  Assume, instead of(1.6), the stronger condition
$$\displaystyle
T-\eta>R_{\Omega}(p).
\tag {1.12}
$$
If $T<2\{(T-\eta)-R_D(p)\}$, then
$$\displaystyle
\lim_{\tau\rightarrow\infty}
e^{\tau T}I_{\partial\Omega}(\tau;B,T)=0.
$$

\endproclaim

Note that the indicator function $I_{\partial\Omega}(\tau;B,T)$ can be computed
from the wave field $u_f$ on $\partial\Omega\times\,]0,\,T[$, generated
by the {\it single} Neumann data $f=f_{B,\,T}$.  Thus, formula (1.11) enables us to know
the quantity $R_D(p)$ which is the radius of the minimum ball centered at $p$ and enclosing $D$.
The point $p$ can be an arbitrary point in $\Bbb R^3$.  We do not mind whether $p\in D$, $p\in\Omega\setminus D$
or $p\in\Bbb R^3\setminus\Omega$.

Condition (1.10) is equivalent to the condition
$$\displaystyle
T>\{(T-\eta)-R_D(p)\}+(R_{\Omega}(p)-R_D(p)).
\tag {1.13}
$$
Under the assumption (1.6) we have
$$\begin{array}{ll}
\displaystyle
2\{(T-\eta)-R_D(p)\}
&
\displaystyle
=\{(T-\eta)-R_D(p)\}+\{(T-\eta)-R_D(p)\}
\\
\\
\displaystyle
&
\displaystyle
\ge
\{(T-\eta)-R_D(p)\}+(R_{\Omega}(p)-R_D(p)).
\end{array}
$$
Therefore, if $T$ satisfies (1.6) and $\displaystyle T>2\{(T-\eta)-R_D(p)\}$,
then $\eta$ satifies (1.13) and hence (1.10).
Thus, assertion (ii) is a direct consequence of (i).

Summing up, we have obtained:
$$\displaystyle
\lim_{\tau\rightarrow\infty}
e^{\tau T}I_{\partial\Omega}(\tau;B,T)
=
\left\{
\begin{array}{ll}
\displaystyle
\infty & \text{if $\eta+R_{\Omega}(p)\le T<2(\eta+R_D(p))$,}
\\
\\
\displaystyle
0 & \text{if $T>2(\eta+R_D(p))$,}
\end{array}
\right.
$$
provided $\eta$ satisfies (1.10).  This criterion gives an alternative and {\it qualitative} characterization of $R_D(p)$ instead of (1.11).

Note that, for all $T$ satisfying (1.6) we have
$$\begin{array}{l}
\,\,\,\,\,\,
\displaystyle
\{(T-\eta)-R_D(p)\}+(R_{\Omega}(p)-R_D(p))\\
\\
\displaystyle
\ge
\inf\,\{\vert P-Q\vert+\vert Q-R\vert\,\vert\,P\in\partial B_{T-\eta}(p),
Q\in\partial D, R\in\partial\Omega\}.
\end{array}
\tag {1.14}
$$
This is proved as follows.  First choose $Q\in\partial D$ such that $R_D(p)=\vert Q-p\vert$.
Second choose $P\in\partial B_{T-\eta}(p)$ such that $Q$ is on the segment $[p,P]$.  Thus we have
$\vert P-Q\vert=(T-\eta)-R_D(p)$.  Third choose $R'\in\partial B_{R_{\Omega}(p)}(p)$ such that $Q$ is on the segment $[p,R']$.
We have $\vert Q-R'\vert=R_{\Omega}(p)-R_{D}(p)$.  Then one can find a point $R\in\partial\Omega$ on the segement $[Q, R']$.
Then we have $\vert Q-R'\vert\ge\vert Q-R\vert$ and thus
$$\displaystyle
\{(T-\eta)-R_D(p)\}+(R_{\Omega}(p)-R_D(p))
=\vert P-Q\vert+\vert Q-R'\vert
\ge \vert P-Q\vert+\vert Q-R\vert.
$$
This yields the desired conclusion.

Note that the right-hand side on (1.14) gives the minimum length of the broken paths that start at $P\in\partial B_{T-\eta}(p)$, 
reflect at $y\in\partial D$ and return to $R\in\partial\Omega$.
Therefore, condition (1.13) is quite natural, and so is (1.10).

If $D$ is large in the sense that $2R_D(p)\ge R_{\Omega}(p)$, then $\eta$ satisfying (1.10) can be arbitrary small.
However, if $2R_D(p)<R_{\Omega}(p)$, then one has to choose a large $\eta$.  The choice of a small $\eta$ depends on a lower estimate
of $R_D(p)$.  This means that we need a-priori information about the size of $R_D(p)$ from below.
However, note that (1.10) is valid for all $\eta$
with $\eta\ge R_{\Omega}(p)$.  This last condition is independent of $D$.

The main difference from \cite{EIV} is the choice of the Neumann data $f$ in (1.1). 
Therein we restrict the location of $B$ to the outside of $\Omega$.
Then the Neumann data in \cite{EIV} is given by
$$\begin{array}{lll}
\displaystyle
f_{B}(x,t)=\frac{\partial}{\partial\nu}v(x,t), & x\in\partial\Omega, & 0\le t\le T,
\end{array}
$$
where $v$ is the solution of (1.2) with this restricted $B$.
So, in this case we have
$$\displaystyle
f_{B,\,T}(x,t)=f_B(x, T-t).
$$
That is, the Neumann data (1.7) plays the role of the {\it time-reversal mirror} \cite{FWCM} equipped on the boundary of $\partial\Omega$
for the wave generated by $f_{B}$ over the time interval $[0,\,T]$ in the case when $D=\emptyset$.
We can generate a natural free wave in $\Omega$ which is emitted on $\partial\Omega$ possibly with some delay, and goes to $B$.

As done in \cite{EIV}, the analysis of the indicator function $I_{\partial\Omega}(\tau;B)$ as $\tau\rightarrow\infty$
is reduced to the study of the asymptotic behaviour of $w^*_{B,T}$ on $D$ as $\tau\rightarrow\infty$.
From (1.2) and (1.9) we know that $w^*_{B,T}$ satisfies
$$\begin{array}{ll}
\displaystyle
(\Delta-\tau^2)w^*_{B,T}-e^{-\tau T}F_0(x)=\partial_tv(x,T)-\tau v(x,T),
&
\displaystyle
x\in\Bbb R^3,
\end{array}
\tag {1.15}
$$
where
$$\begin{array}{ll}
\displaystyle
F_0(x)=-\Psi_B(x), & x\in\Bbb R^3.
\end{array}
\tag {1.16}
$$
Changing the role of $F_0$ and $\partial_tv(x,T)-\tau v(x,T)$ in (1.15), we have
$$\begin{array}{ll}
\displaystyle
(\Delta-\tau^2)w^*_{B,T}+(\tau v(x,T)-\partial_tv(x,T))=e^{-\tau T}F_0,
&
\displaystyle
x\in\Bbb R^3.
\end{array}
$$
Since the trem $e^{-\tau T}F_0$ can be ignored
as $\tau\rightarrow\infty$, we have to consider
the term $\tau v(x,T)-\partial_tv(x,T)$, the main source.
This together with (1.4) and (1.5) leads us to study the asymptotic
profile of the following integral as $\tau\rightarrow\infty$:
$$\begin{array}{ll}
\displaystyle
\frac{1}{4\pi}
\int_{B_{T+\eta}(p)\setminus B_{T-\eta}(p)}
\frac{e^{-\tau\vert x-y\vert}}{\vert x-y\vert}
(\tau v(y,T)-\partial_tv(y,T))\,dy, & x\in D.
\end{array}
$$
This is a new situation not appearing in \cite{EIV}.
In this paper, using Kirchhoff's formula
(1.3), we compute the integral explicitly and clarify the asymptotic behaviour
as $\tau\rightarrow\infty$.

The procedure for extracting $R_D(p)$ is explicit, direct and has the following feature:
in the processing of the signal we do not make use of the knowledge of the boundary condition.
Note that in contrast to this,
the so-called continuation procedure of the solutions of the governing equation close to obstacle
makes use of the boundary condition of the obstacle in the procedure, 
such as that of \cite{LCW} and also \cite{BP}, which is a combination
of a continuation method in the {\it frequency domain} and the Fourier transform.

A numerical method in \cite{dBK} for a {\it penetrable obstacle} (embedded in the whole plane $\Bbb R^2$) is a combination
of a time-reversed scattered wave field continuation method and an optimization method for unknown wave speed in the obstacle.
To continue the scattered wave field from the obstacle, they choose a disc that {\it encloses} the obstacle
and solve numerically a time-reversed initial boundary value problem
for the original governing equation in an annulus like domain whose inner boundary
is the boundary of the disc with a time reversed absorbing boundary condition.  On the outer boundary
of the domain where the observed data are collected, the time-reversed scattered field is prescribed as another boundary condition.
Using the computed scattered field in the annulus domain, they introduce an optimization problem with respect to the unknown wave speed
in the obstacle.  It seems that it is not clear whether their method can cover the case when the wave only propagate in a bounded domain,
not the whole space like our situation, since in that case one has to consider the scattered wave not only from the obstacle
but also from the outer boundary.

We mention an analytical approach due to Oksanen \cite{O} which is based on the {\it boundary control method}, see \cite{BEL, B2}.
Therein a similar inverse obstacle problem for the wave equation in a bounded domain or compact manifold with
a boundary is considered.
The approach therein enables us to compute the volume of a set, called the {\it domain of influence},
which is closely related to an unknown obstacle embedded in the domain.
Intuitively, in our Euclidean setting,
it is the set of all points $x\in\Omega\setminus\overline D$ such that the wave, governed by the wave equation in $\Omega\setminus\overline D$
generated at some point $y_0$ on $\partial\Omega$ at $t=0$,
reaches at $x$ within the time $T(y_0)$, where $T(y)$, $y\in\partial\Omega$ is an arbitrary given continuous 
function with the values in $[0,\,T/2]$
and $T(y)=0$ for $y\in\partial\Omega\setminus\overline \Gamma$; $\Gamma$ is an arbitrary prescribed non empty open subset of $\partial\Omega$.
The computed volume yields some information about the location of the obstacle.
The point is to construct a one parameter family of the Neumann data $f$ in such a way that $u_f(x,T/2)$
{\it approximates} the characteristic function of the domain of influence.
The construction is reduced to solving an equation with a parameter written by the Tikhonov regularization of a linear operator
on the boundary of the domain.
The operator is written by using the local hyperbolic Neumann-to-Dirichlet operator and {\it time-reversal operation} on the boundary.
It appears in Blagovestchenskii's identity and is the base of the boundary control method (cf. \cite{B2}).
The idea of the construction is closely related to the {\it focusing wave approach} developed for the wave speed determination problem,
see \cite{BKLS, DKL} and references therein.
However, his result does not tell us what information about the unknown obstacle can be extracted from a {\it single set} of the Dirichlet and Neumann data.
Note that in the crucial step of the proof for the justification of his method, the unique continuation property of 
the governing equation is essential even in our simple situation.  
Our method together with the proof is free from the property, simple and rather elementary.

The {\it linear sampling method} in the time domain has been developed, for example,
in \cite{MS} for an inverse obstacle problem
in waveguide geometry.  The method employs
output data corresponding to infinitely many input.
The output data for each input is
observed taken over the infinite time interval $0<t<\infty$.
Thus, the time-reversal operation never appears.

Finally, we point out that, in \cite{IK2} an extraction formula of $R_D(p)$ is given when $D$ is an inclusion 
embedded in a homogeneous isotropic conductive medium, and the governing equation of the signal propagating inside the medium
is given by the {\it heat equation}.  
It is easy to see that the result therein also covers the cavity case treated in \cite{IK1}.
The data used therein is the Neumann-to-Dirichlet map in the time domain acting on the special Neumann data
having the separation of variables form
$$\displaystyle
\varphi(t)\,\frac{\partial v_{\tau}}{\partial\nu}(x;p),
$$
where $p$ is an arbitrary point in $\Bbb R^3$, say $\varphi(t)\sim t^m$ as $t\downarrow 0$ for an integer $m$, and 
$$
\displaystyle
v_{\tau}(x;p)=
\left\{
\begin{array}{ll}
\displaystyle
\frac{\sinh\sqrt{\tau}\vert x-p\vert}{\vert x-p\vert}, & x\in\Bbb R^3\setminus\{p\},\\
\\
\displaystyle
\sqrt{\tau}, & x=p.
\end{array}
\right.
\tag {1.17}
$$
Since $v_{\tau}$ depends on $\tau$, in this sense, the data to determine $R_D(p)$ for a fixed $p$ is {\it infinitely many}.
In this sense the result shares the same spirit
as a typical result in the classical enclosure method \cite{E00}, which employs infinitely many observation data.
However, note that the normal derivative of $v_{\tau}$ blows up as $\tau\rightarrow\infty$.
It should be emphasized that the Neumann data $f_{B,T}$ given by (1.7) is {\it independent} of such a parameter
which causes the blowing up.
At the present time the author does not know whether there exists a suitable Neumann data depending {\it only} on $p$ or a ball
centered at $p$ with a small radius that
yields $R_D(p)$ for inverse obstacle problems governed by the heat equation.
The main obstruction is the lack of time-reversal invariance and that of the occurrence of lacuna for the fundamental solution.

A brief outline of this paper is as follows.  Theorem 1.1 is proved in Section 2.  
The proof starts with describing the decomposition formula of the indication function.
Using this formula, together with a lemma concerning an upper bound for the second term in the formula,
we reduce the problem to deriving estimates of the energy integral for $w_{B,T}^*$ 
as $\tau\rightarrow\infty$ from above and below.
For the purpose, using the time domain expression (1.3) of $v$, we explicitly write the leading profile of
$w_{B,T}^*$ in $B_{T-\eta}(p)$ as $\tau\rightarrow\infty$ down as stated in Lemma 2.2.
This is the key point of the proof of Theorem 1.1.  
The proof of Lemma 2.2 is given in Section 3.  Since the proof requires explicit forms of some volume integrals,
we give their derivation in Appendix.

\section{Proof of Theorem 1.1.}

In this section, for simlicity of description, we always write
$$\begin{array}{lll}
\displaystyle
w=w_{B,T},
& 
\displaystyle
w^*=w_{B,T}^*,
&
\displaystyle
R=w-w^*.
\end{array}
$$
The following decomposition formula is valid
(see Proposition 2.1 in \cite{EIV}). 

\proclaim{\noindent Proposition 2.1.}
We have
$$\displaystyle
I_{\partial\Omega}(\tau;B,T)
=J_*(\tau)+E(\tau)+{\cal R}(\tau),
\tag {2.1}
$$
where
$$\displaystyle
J_*(\tau)=\int_D(\vert\nabla w^*\vert^2+\tau^2\vert w^*\vert^2)\,dx,
\tag {2.2}
$$
$$\displaystyle
E(\tau)
=\int_{\Omega\setminus\overline D}
(\vert\nabla R\vert^2+\tau^2\vert R\vert^2)\,dx,
\tag {2.3}
$$
$$\displaystyle
{\cal R}(\tau)
=e^{-\tau T}
\left\{
\int_DF_0w^*dx+\int_{\Omega\setminus\overline{D}}FRdx+\int_{\Omega\setminus\overline D}(F_0-F)w^*dx\right\},
$$
$$\begin{array}{ll}
\displaystyle
F=F(x,\tau)=\partial_tu_f(x,T)+\tau u_f(x,T),
& 
\displaystyle
x\in\Omega\setminus\overline D
\end{array}
\tag {2.4}
$$
and $F_0$ is given by (1.16).

\endproclaim

Note that the proof of (2.1) is based on the two facts.

First it follows from (1.1) and (1.8) that $w$ satisfies
$$\left\{
\begin{array}{ll}
\displaystyle
(\Delta-\tau^2)w=e^{-\tau T}F & \text{in $\Omega\setminus\overline D$,}
\\
\\
\displaystyle
\frac{\partial w}{\partial\nu}=\frac{\partial w^*}{\partial\nu}
& \text{on $\partial\Omega$},
\\
\\
\displaystyle
\frac{\partial w}{\partial\nu}=0 &
\text{on $\partial D$.}
\end{array}
\right.
\tag {2.5}
$$
This is the same as before.
Then, from (1.5), assumption (1.6) and (1.15) we see that $w^*$ satisfies
$$\begin{array}{ll}
\displaystyle
(\Delta-\tau^2)w^*=e^{-\tau T} F_0, & x\in\Omega.
\end{array}
\tag {2.6}
$$
Using (2.5) and (2.6) together with integration by parts we obtain (2.1).

Similar to Lemma 2.2 in \cite{EIV}, we have

\proclaim{\noindent Lemma 2.1(Dominance estimate).}
We have
$$\displaystyle
E(\tau)
=O(\tau^2 J_*(\tau)+\tau^2 e^{-2\tau T})
\tag {2.7}
$$
as $\tau\rightarrow\infty$.
\endproclaim

The point is: $f=f_{B,T}$ is \underline{independent} of $\tau$ and thus (2.4) gives $\Vert F\Vert_{L^2(\Omega\setminus\overline D)}=O(\tau)$.
This, together with (1.16) gives  $\Vert F-F_0\Vert_{L^2(\Omega\setminus\overline D)}=O(\tau)$, which is the
same as in the proof of Lemma 2.2 in \cite{EIV}.

The next task is to give an upper bound on ${\cal R}(\tau)$ and the upper and lower estimates on $J_*(\tau)$.
For the purpose, we study the local behaviour of $w^*$.

From (1.15) we see that the $w^*$ takes the form
$$\displaystyle
w^*=w_1^*+e^{-\tau T}w_R^*,
\tag {2.8}
$$
where
$$\begin{array}{ll}
\displaystyle
w_1^*(x,\tau)=\frac{1}{4\pi}\int_{\Bbb R^3}\frac{e^{-\tau\vert x-y\vert}}{\vert x-y\vert}\,(\tau v(y,T)-\partial_t v(y,T))dy,
& x\in\Bbb R^3
\end{array}
\tag {2.9}
$$
and $w_R^*$ satisfies
$$\begin{array}{ll}
\displaystyle
(\Delta-\tau^2)w_R^*+\Psi_B=0, & x\in\Bbb R^3.
\end{array}
$$
By integration by parts we have immediately, as $\tau\rightarrow\infty$,
$$
\displaystyle
\tau\Vert w_R^*\Vert_{L^2(\Bbb R^3)}+\Vert\nabla w_R^*\Vert_{L^2(\Bbb R^3)}
=O(1).
\tag {2.10}
$$

Thus, to clarify the behaviour of $w^*$ in $B_{T-\eta}(p)$ it suffices to study
that of $w_1^*$.  Noting (1.4) and (1.5), we prepare
two lemmas in which the first one yields an explicit form of 
$w_1^*(x,\tau)$ for $x\in B_{T-\eta}(p)$ and the second its upper and lower estimates.

\proclaim{\noindent Lemma 2.2.}
Let $T>\eta$.
We have
$$\begin{array}{l}
\,\,\,\,\,\,\displaystyle
\frac{\tau^2}{4\pi}
\int_{B_{T+\eta}(p)\setminus B_{T-\eta}(p)}
\frac{e^{-\tau\vert x-y\vert}}{\vert x-y\vert}
(\tau v(y,T)-\partial_tv(y,T))\,dy\\
\\
\displaystyle
=
e^{-\tau(T-\eta)}{\cal H}(\tau;T,\eta)
\,\frac{\sinh\tau\vert x-p\vert}{\vert x-p\vert}
\end{array}
$$
for all $x\in B_{T-\eta}(p)\setminus\{p\}$,
where
$$\displaystyle
{\cal H}(\tau;T,\eta)
=\tau^{-1}\left(\eta+O(\tau^{-1})\right).
$$

\endproclaim

For the proof of Lemma 2.2 see Section 3.  It is a chain of a careful explicit computation
using the speciality of the form.

Note that the function
$$\begin{array}{ll}\displaystyle
\frac{\sinh\tau\vert x-p\vert}{\vert x-p\vert}, & x\in\Bbb R^3\setminus\{p\}
\end{array}
$$
has a unique extension to the whole space as a smooth function and satisfies
the modified Helmholtz equation $(\Delta-\tau^2)v=0$ in the whole space. 
More precisely, the function coincideds with $v_{\tau^2}(x;p)$, which is given by (1.17), with $\tau$ replaced with $\tau^2$. 
In the following lemma we continue to use this notation to denote its extension.

\proclaim{\noindent Lemma 2.3.}
Let $U$ be an arbitrary bounded open subset of $\Bbb R^3$
and $p$ an arbitrary point in $\Bbb R^3$.  Set $R_U(p)=\sup_{x\in U}\,\vert x-p\vert$.

(i)  There exists a real number $\mu_1$ such that, as $\tau\rightarrow\infty$
$$\displaystyle
\int_U v_{\tau^2}(x;p)^2\,dx
+\int_U\left\vert\nabla v_{\tau^2}(x;p)\right\vert^2\,dx
=O(\tau^{2\mu_1}e^{2\tau R_U(p)}).
$$

(ii)  Assume that $\partial U$ is Lipschitz.
There exist positive numbers $C$ and $\tau_0$ and a real number $\mu_2$ such that
$$\displaystyle
\tau^{2\mu_2}e^{-2\tau R_U(p)}\,\int_U v_{\tau^2}(x;p)^2\,dx
\ge C
$$
for all $\tau\ge\tau_0$.

\endproclaim

Since $U\subset B_{R_U(p)}(p)$, the proof of Lemma 2.3 (i) can be done by replacing $U$ with the ball $B_{R_U(p)}(p)$ 
and using the polar coordinates around $p$.
The proof of Lemma 2.3 (ii) can be done by using the same argument for the proof of Lemma 6 in \cite{II}.
The point of the argument is to find a subdomain $\tilde{U}$ of $U$ such that $R_{\tilde{U}}(p)=R_U(p)$ and
$\vert x-p\vert\ge {R_U(p)}/2$ for all $x\in\tilde{U}$.  For this purpose, the Lipschitz regularity of $\partial U$
is enough.  For these reasons, we omit the proof of Lemma 2.3.
Note that the concrete values of $\mu_1$ and $\mu_2$ are not essential in this paper just like in \cite{II}
and other papers for the time domain enclosure methods.

From (1.4), (1.5), expression (2.9) and Lemma 2.2, one gets an explicit asymptotic form of $w_1^*$ in $B_{T-\eta}(p)$.
Then, from Lemma 2.3 together with (2.8) and (2.10), we immediately obtain

\proclaim{\noindent Lemma 2.4(Propagation estimate).}
Let $U$ be an arbitrary bounded open subset of $\Bbb R^3$ such that
$U\subset B_{T-\eta}(p)$, that is,
$$\displaystyle
T-\eta\ge R_U(p).
\tag {2.11}
$$

(i)  There exist a real number $\mu_3$ such that, as $\tau\rightarrow\infty$
$$\displaystyle
\tau\Vert w^*\Vert_{L^2(U)}+\Vert\nabla w^*\Vert_{L^2(U)}=O(\tau^{\mu_3}e^{-\tau(T-\eta)}e^{\tau R_{U}(p)}+\tau^2 e^{-\tau T}).
$$

(ii) If $\partial U$ is Lipschitz, then there exist positive numbers $\tau_0$ and $C$ such that
$$\displaystyle
\tau^{\mu_2+3}e^{\tau(T-\eta)}e^{-\tau R_U(p)}\Vert w^*\Vert_{L^2(U)}\ge C
$$
for all $\tau\ge\tau_0$, where $\mu_2$ is the same as that of Lemma 2.3 (ii).
\endproclaim

Using the facts $\Vert F\Vert_{L^2(\Omega\setminus\overline D)}=O(\tau)$ and $\Vert F_0\Vert_{L^2(D)}=O(1)$
together with (2.2), (2.3) and (2.7) we have, as $\tau\rightarrow\infty$
$$\left\{
\begin{array}{l}
\displaystyle
\int_{\Omega\setminus\overline D}FRdx
=O(\tau\cdot\tau^{-1}E(\tau)^{1/2})=O(E(\tau)^{1/2})=O(\tau J_*(\tau)^{1/2}+\tau e^{-\tau T}),
\\
\\
\displaystyle
\int_DF_0w^*dx
=O(\tau^{-1}J_*(\tau)^{1/2}).
\end{array}
\right.
$$
Applying Lemma 2.4 (i) to the case $U=\Omega$, we obtain
$$\displaystyle
\Vert w^*\Vert_{L^2(\Omega\setminus\overline D)}
=O(\tau^{\mu_3-1}e^{-\tau(T-\eta)}e^{\tau R_{\Omega}(p)}
+\tau e^{-\tau T})
$$
and this, together with $\Vert F-F_0\Vert_{L^2(\Omega\setminus\overline D)}
=O(\tau)$, yields
$$\displaystyle
\int_{\Omega\setminus\overline D}(F_0-F)w^*\,dx
=O(\tau^{\mu_3}e^{-\tau(T-\eta)}e^{\tau R_{\Omega}(p)}+\tau^2 e^{-\tau T}).
$$
Moreover, from Lemma 2.4 (i) in the case $U=D$, we obtain
$$\displaystyle
J_*(\tau)=O(\tau^{2\mu_3} e^{-2\tau(T-\eta)}e^{2\tau R_D(p)}+\tau^4 e^{-2\tau T}).
\tag {2.12}
$$
From these, we obtain
$$\begin{array}{l}
\displaystyle
\,\,\,\,\,\,
{\cal R}(\tau)\\
\\
\displaystyle
=O(e^{-\tau T}
(\tau J_*(\tau)^{1/2}+\tau e^{-\tau T}))
+O(\tau^{\mu_3}e^{-\tau T}e^{-\tau(T-\eta)}e^{\tau R_{\Omega}(p)}+\tau^2 e^{-2\tau T}))\\
\\
\displaystyle
=O(e^{-\tau T}
\left\{\tau(\tau^{\mu_3} e^{-\tau(T-\eta)}e^{\tau R_D(p)}+\tau^2e^{-\tau T})+\tau e^{-\tau T}\right\})\\
\\
\displaystyle
\,\,\,
+O(\tau^{\mu_3}e^{-\tau T}e^{-\tau(T-\eta)}e^{\tau R_{\Omega}(p)}+\tau^2 e^{-2\tau T})\\
\\
\displaystyle
=O(\tau^{\mu_3+1} e^{-\tau T}e^{-\tau(T-\eta)}e^{\tau R_D(p)}+\tau^{2}e^{-2\tau T}+\tau^{\mu_3}e^{-\tau T}e^{-\tau(T-\eta)}e^{\tau R_{\Omega}(p)}).
\end{array}
\tag {2.13}
$$
Thus
$$\begin{array}{l}
\displaystyle
\,\,\,\,\,\,
e^{2\tau(T-\eta)}e^{-2\tau R_D(p)}{\cal R}(\tau)
\\
\\
\displaystyle
=e^{2\tau(T-\eta)}e^{-2\tau R_D(p)}O(\tau^{\mu_3+1}e^{-\tau T}e^{-\tau(T-\eta)}e^{\tau R_D(p)}
+\tau^{2}e^{-2\tau T}+\tau^{\mu_3}e^{-\tau T}e^{-\tau(T-\eta)}e^{\tau R_{\Omega}(p)})\\
\\
\displaystyle
=O(\tau^{\mu_3+1} e^{-\tau\eta}e^{-\tau R_D(p)}+\tau^{2}e^{-2\tau\eta}e^{-2\tau R_D(p)}+\tau^{\mu_3}e^{-\tau(\eta+2R_D(p)-R_{\Omega}(p))}).
\end{array}
\tag {2.14}
$$

Now we are ready to describe the proof of Theorem 1.1 (i).
Let $\eta$ satisfy the condition (1.10).
Then (2.14) yields
$$\displaystyle
e^{2\tau(T-\eta)}e^{-2\tau R_D(p)}{\cal R}(\tau)=O(\tau^{-\infty}).
\tag {2.15}
$$
Thus, from this, (2.1), (2.7) and (2.12), we obtain
$$\begin{array}{l}
\displaystyle
\,\,\,\,\,\,
I_{\partial\Omega}(\tau;B,T)\\
\\
\displaystyle
=O(\tau^{2+2\mu_3}e^{-2\tau(T-\eta)}e^{2\tau R_D(p)}+\tau^6 e^{-2\tau T})
\\
\\
\displaystyle
=O(\tau^{2+2\mu_3} e^{-2\tau(T-\eta)}e^{2\tau R_D(p)}
(1+\tau^{6-2-2\mu_3}e^{-\tau(\eta+R_D(p))}))
\\
\\
\displaystyle
=O(\tau^{2+2\mu_3} e^{-2\tau(T-\eta)}e^{2\tau R_D(p)}).
\end{array}
\tag {2.16}
$$
Moreover, from (2.1), (2.2) and (2.15) we have
$$\begin{array}{ll}
\displaystyle
e^{2\tau(T-\eta)}e^{-2\tau R_D(p)}I_{\partial\Omega}(\tau;B,T)
&
\displaystyle
\ge e^{2\tau(T-\eta)}e^{-2\tau R_D(p)}J_*(\tau)+O(\tau^{-\infty})\\
\\
\displaystyle
&
\displaystyle
\ge \tau^2 e^{2\tau(T-\eta)}e^{-2\tau R_D(p)}\Vert w^*\Vert_{L^2(D)}^2+O(\tau^{-\infty}).
\end{array}
$$
Since (1.6) implies (2.11) with $U=D$, from Lemma 2.4 (ii) in the case when $U=D$
and writing
$$\displaystyle
\tau^2=\tau^{-2(\mu_2+2)}\tau^{2(\mu_2+3)},
$$
one can conclude that
there exist positive numbers $C$ and $\tau_0$ such that
$$
\displaystyle
\tau^{2(\mu_2+2)}e^{2\tau(T-\eta)}e^{-2\tau R_D(p)}I_{\partial\Omega}(\tau;B,T)
\ge C
\tag {2.17}
$$
for all $\tau\ge\tau_0$.
A combination of (2.16) and (2.17) ensures that assertion (i) is valid.

As pointed out in the second paragraph following Theorem 1.1, assertion (ii) is a direct consequence of (i).
Thus, it suffices to prove (iii).
Instead of (2.15) which is a consequence of assumption (1.10),
we go back to (2.13).
Then we have
$$\begin{array}{ll}
\displaystyle
e^{\tau T}{\cal R}(\tau)
&
\displaystyle
=O(\tau^{\mu_3+1} e^{-\tau(T-\eta)}e^{\tau R_D(p)}+\tau^{2}e^{-\tau T}+\tau^{\mu_3}e^{-\tau(T-\eta)}e^{\tau R_{\Omega}(p)}).
\end{array}
$$
Note that $T<2\{(T-\eta)-R_D(p)\}$ implies that $T-\eta>R_D(p)$.
Thus under assumption (1.12), which is stronger than (1.6), we conclude
$$\displaystyle
e^{\tau T}{\cal R}(\tau)=O(\tau^{-\infty}).
$$
Now from this, (2.1), (2.7) and (2.12) we obtain
$$
\displaystyle
e^{\tau T}I_{\partial\Omega}(\tau;B,T)
=O(\tau^{2\mu_3+2}e^{\tau T} e^{-2\tau(T-\eta)}e^{2\tau R_D(p)}\,)+O(\tau^{-\infty}).
$$
Since $T<2\{(T-\eta)-R_D(p)\}$, we conclude
$$\displaystyle
e^{\tau T}I_{\partial\Omega}(\tau;B,T)=O(\tau^{-\infty}).
$$

This completes the proof of Theorem 1.1.

\section{Proof of Lemma 2.2}

First we compute the value of $v(x,T)$, together with $\partial_t v(x,T)$,
at $x\in B_{T+\eta}(p)\setminus B_{T-\eta}(p)$.

\proclaim{\noindent Proposition 3.1.}

(i)  If $\vert\vert x-p\vert-t\vert<\eta$ and $\eta<\vert x-p\vert+t$, then we have
$$
\left\{
\begin{array}{l}
\displaystyle
v(x,t)
=\frac{1}{2}
\left\{
\frac{\eta^3}{6\vert x-p\vert}
-
\frac{\eta(\vert x-p\vert-t)^2}{2\vert x-p\vert}
+\frac{\vert \vert x-p\vert-t\vert^3}{3\vert x-p\vert}
\right\},
\\
\\
\displaystyle
\partial_tv(x,t)
=\frac{\vert x-p\vert-t}{2\vert x-p\vert}
\left(\eta-\vert\vert x-p\vert-t\vert\right).
\end{array}
\right.
$$

(ii) If $\vert x-p\vert+t<\eta$, then
we have
$$\left\{
\begin{array}{l}
\displaystyle
v(x,t)
=\eta t
-\frac{1}{6\vert x-p\vert}
\left\{(\vert x-p\vert+t)^3-\vert\vert x-p\vert-t\vert^3\right\},
\\
\\
\displaystyle
\partial_tv(x,t)
=\eta
-\frac{1}{2\vert x-p\vert}
\left\{(\vert x-p\vert+t)^2+(\vert x-p\vert-t)\vert\vert x-p\vert-t\vert\right\}.
\end{array}
\right.
$$

\endproclaim

{\it\noindent Proof.}
Write (1.3) as
$$\displaystyle
v(x,t)=\frac{t}{4\pi}
\int_{S(x;B)}(\eta-\vert(x+t\omega)-p\vert)\,d\omega,
\tag {3.1}
$$
where
$$\displaystyle
S(x;B)=\{
\omega\in S^2\,\vert\,
\vert (x+t\omega)-p\vert<\eta
\}.
$$
The inequality $\vert (x+t\omega)-p\vert<\eta$ for $\omega\in S^2$ is equivalent to
$$\displaystyle
\omega\cdot\frac{p-x}{\vert p-x\vert}
>\frac{\vert p-x\vert^2+t^2-\eta^2}{2t\vert p-x\vert}.
$$

First consider the case when $\vert \vert x-p\vert-t\vert<\eta$ and $\vert x-p\vert+t>\eta$.
In this case, we have
$$\displaystyle
-1<\frac{\vert p-x\vert^2+t^2-\eta^2}{2t\vert p-x\vert}<1.
$$
Define
$$\displaystyle
\phi_0=\arccos\frac{\vert p-x\vert^2+t^2-\eta^2}{2t\vert p-x\vert}.
$$
Then one can write all the points $\omega\in S(x;B)$ in terms of the polar coordinates:
$$\displaystyle
\omega=\sin\phi\,(\cos\theta\,\mbox{\boldmath $b$}+\sin\theta\,\mbox{\boldmath $c$})+\cos\phi\,\frac{p-x}{\vert p-x\vert},
$$
where $0\le\theta\le 2\pi$, $0\le\phi<\phi_0$ and the
unit vectors $\mbox{\boldmath $b$}$ and $\mbox{\boldmath $c$}$ are parpendicular
and satisfy 
$$\displaystyle
\mbox{\boldmath $b$}\times\mbox{\boldmath $c$}=\frac{p-x}{\vert p-x\vert}.
$$
Thus, (3.1) becomes
$$\displaystyle
v(x,t)
=\frac{t}{2}
\int_0^{\phi_0}\sin\phi
\,(\eta-\sqrt{\vert x-p\vert^2+t^2-2t\vert x-p\vert\cos\phi})\,d\phi.
\tag {3.2}
$$
Here we have
$$\begin{array}{ll}
\displaystyle
\int_0^{\phi_0}\sin\phi\,d\phi
&
\displaystyle
=1-\cos\phi_0\\
\\
\displaystyle
&
\displaystyle
=1-\frac{\vert x-p\vert^2+t^2-\eta^2}{2t\vert x-p\vert}
\\
\\
\displaystyle
&
\displaystyle
=\frac{\eta^2-(\vert x-p\vert-t)^2}{2t\vert x-p\vert}
\end{array}
$$
and
$$
\begin{array}{l}
\displaystyle
\,\,\,\,\,\,
\int_0^{\phi_0}\sqrt{\vert x-p\vert^2+t^2-2t\vert x-p\vert\cos\phi}\,\sin\phi\,d\phi\\
\\
\displaystyle
=\frac{1}{3t\vert x-p\vert}
(\sqrt{\vert x-p\vert^2+t^2-2t\vert x-p\vert\,\cos\phi}\,)^3\vert_{\phi=0}^{\phi=\phi_0}\\
\\
\displaystyle
=\frac{1}{3t\vert x-p\vert}
\left\{
(\sqrt{\vert x-p\vert^2+t^2-2t\vert x-p\vert\,\cos\phi_0}\,)^3
-\vert \vert x-p\vert-t\vert^3)
\right\}\\
\\
\displaystyle
=\frac{1}{3t\vert x-p\vert}
\left(\eta^3-\vert\vert x-p\vert-t\vert^3\right).
\end{array}
$$
Thus, from (3.2), we have
$$
\displaystyle
v(x,t)
=\frac{t}{2}\cdot\eta\cdot\frac{\eta^2-(\vert x-p\vert-t)^2}{2t\vert x-p\vert}
-\frac{t}{2}\cdot
\frac{1}{3t\vert x-p\vert}
\left(\eta^3-\vert\vert x-p\vert-t\vert^3\right).
$$
We have
$$
\displaystyle
\partial_t(\vert\vert x-p\vert-t\vert^3)=-3(\vert x-p\vert-t)\vert\vert x-p\vert-t\vert.
\tag {3.3}
$$
Thus, one gets
$$
\displaystyle
\partial_tv(x,t)
=
\frac{\eta(\vert x-p\vert-t)}{2\vert x-p\vert}
-\frac{(\vert x-p\vert-t)\vert\vert x-p\vert-t\vert
}{2\vert x-p\vert}.
$$
This yields the desired conclusion of (i).

Next consider the case when $\vert x-p\vert+t<\eta$.  We see that
$$\displaystyle
\frac{\vert p-x\vert^2+t^2-\eta^2}{2t\vert p-x\vert}<-1.
$$
Thus, $S(x;B)=S^2$ and, using the same polar coordinates as above with $\phi_0=\pi$, we have
$$\displaystyle
v(x,t)=\frac{t}{2}\cdot2\eta
-\frac{t}{2}\cdot \frac{1}{3t\vert x-p\vert}
\left\{(\vert x-p\vert+t)^3-\vert\vert x-p\vert-t\vert^3\right\}.
$$
Using (3.3), we have
$$\displaystyle
\partial_tv(x,t)
=\eta-\frac{1}{6\vert x-p\vert}
\left\{
3(\vert x-p\vert+t)^2+
3(\vert x-p\vert-t)\vert\vert x-p\vert-t\vert\right\}.
$$
This yields the desired formula in (ii).

\noindent
$\Box$

{\bf\noindent Remark 3.1.}
(a)  Let $\vert \vert x-p\vert-t\vert<\eta$ and $\eta<\vert x-p\vert+t$.
Then, we have
$$\displaystyle
\vert (\vert x-p\vert-t)-\eta\vert<2t.
$$
Thus, from Proposition 3.1 (i), we have $\lim_{t\downarrow 0}\partial_tv(x,t)=0$.

(b)  Let $\vert x-p\vert<\eta$. Then, for all $t>0$ with
$t<\eta-\vert x-p\vert$, we have $\vert x-p\vert+t<\eta$.  Then form Proposition 3.1 (ii) we obtain
$\lim_{t\downarrow 0}\partial_tv(x,t)=\eta-\vert x-p\vert$.

Thus, it suffices to compute the integrals
$$\begin{array}{ll}
\displaystyle
I_j(x;R_1,R_2)
=\frac{1}{4\pi}
\int_{B_{R_2}(p)\setminus B_{R_1}(p)}
\frac{e^{-\tau\vert x-y\vert}}{\vert x-y\vert}
\,\vert y-p\vert^j\,dy, & x\in B_{R_1}(p)
\end{array}
$$
for $j=-1, 0, 1, 2$ and $R_2>R_1$.

The resluts are listed below which are the direct consequence of Proposition A in Appendix:
$$\begin{array}{ll}
\displaystyle
I_j(x;R_1,R_2)
=\frac{1}{\tau^2}H_j(\tau;R_1,R_2)\frac{\sinh\tau\vert x-p\vert}{\vert x-p\vert},
&
x\in B_{R_1}(p)\setminus\{p\},
\end{array}
\tag {3.4}
$$
where
$$
\left\{
\begin{array}{l}
\displaystyle
H_{-1}(\tau;R_1,R_2)
=e^{-\tau R_1}-e^{-\tau R_2},
\\
\\
\displaystyle
H_{0}(\tau;R_1,R_2)
=\left(R_1+\frac{1}{\tau}\right)e^{-\tau R_1}
-\left(R_2+\frac{1}{\tau}\right)e^{-\tau R_2},
\\
\\
\displaystyle
H_{1}(\tau;R_1,R_2)
=\left(R_1^2+\frac{2}{\tau}
R_1+\frac{2}{\tau^2}\right)e^{-\tau R_1}
-\left(R_2^2+\frac{2}{\tau}
R_2+\frac{2}{\tau^2}\right)e^{-\tau R_2},
\\
\\
\displaystyle
H_{2}(\tau;R_1,R_2)
=\left(R_1^3+\frac{3}{\tau}R_1^2+\frac{6}{\tau^2}R_1+\frac{6}{\tau^3}\right)e^{-\tau R_1}
-
\left(R_2^3+\frac{3}{\tau}R_2^2+\frac{6}{\tau^2}R_2+\frac{6}{\tau^3}\right)e^{-\tau R_2}.
\end{array}
\right.
$$

From Proposition 3.1 and (3.4) we obtain the following result.

\proclaim{\noindent Proposition 3.2.}
Let $T>\eta$. 
We have the expression
$$\begin{array}{l}
\,\,\,\,\,\,\displaystyle
\frac{1}{4\pi}
\int_{B_{T+\eta}(p)\setminus B_{T-\eta}(p)}
\frac{e^{-\tau\vert x-y\vert}}{\vert x-y\vert}
(\tau v(y,T)-\partial_tv(y,T))\,dy\\
\\
\displaystyle
=\frac{1}{\tau^2}({\cal H}_{+}(\tau;T,\eta)+{\cal H}_{-}(\tau;T,\eta))\,\frac{\sinh\tau\vert x-p\vert}{\vert x-p\vert}
\end{array}
$$
for all $x\in B_{T-\eta}(p)\setminus\{p\}$, where
$$\begin{array}{l}
\displaystyle
\,\,\,\,\,\,
{\cal H}_{+}(\tau;T,\eta)\\
\\
\displaystyle
=
\left\{\frac{1}{12}\tau(\eta-2T)(\eta+T)^2
+\frac{1}{2}T(\eta+T)\right\}H_{-1}(\tau;T,T+\eta)\\
\\
\displaystyle
\,\,\,
+\left\{\frac{1}{2}\tau T(\eta+T)-\frac{1}{2}(\eta+2T)\right\}H_{0}(\tau;T,T+\eta)\\
\\
\displaystyle
\,\,\,
+\left\{-\frac{1}{4}\tau(\eta+2T)+\frac{1}{2}\right\}H_{1}(\tau;T,T+\eta)
+\frac{1}{6}\tau H_{2}(\tau;T,T+\eta)
\end{array}
\tag{3.5}
$$
and
$$\begin{array}{l}
\displaystyle
\,\,\,\,\,\,
{\cal H}_{-}(\tau;T,\eta)
\\
\\
\displaystyle
=
\left\{\frac{1}{12}\tau(\eta+2T)(\eta-T)^2
+\frac{1}{2}T(\eta-T)\right\}H_{-1}(\tau;T-\eta,T)\\
\\
\displaystyle
\,\,\,
+\left\{\frac{1}{2}\tau T(\eta-T)-\frac{1}{2}(\eta-2T)\right\}
H_0(\tau;T-\eta,T)\\
\\
\displaystyle
\,\,\,
+\left\{-\frac{1}{4}\tau(\eta-2T)-\frac{1}{2}\right\}H_1(\tau;T-\eta,T)
-\frac{1}{6}\tau H_2(\tau;T-\eta,T).
\end{array}
\tag{3.6}
$$

\endproclaim

{\it\noindent Proof.}
Consider the decomposition
$$\displaystyle
B_{T+\eta}(p)\setminus B_{T-\eta}(p)=B_1\cup B_2,
$$
where
$$\displaystyle
B_1=\{y\in\Bbb R^3\,\vert\,T\le\vert y-p\vert<T+\eta\}=B_{T+\eta}(p)\setminus B_T(p)
$$
and
$$\displaystyle
B_2=\{y\in\Bbb R^3\,\vert\,T-\eta<\vert y-p\vert\le T\}=B_T(p)\setminus B_{T-\eta}(p).
$$
Since $T>\eta$, Proposition 3.1 yields
$$
\displaystyle
v(y,T)
=
\left\{
\begin{array}{ll}
\displaystyle
\frac{1}{2}
\left\{
\frac{\eta^3}{6\vert y-p\vert}
-
\frac{\eta(\vert y-p\vert-T)^2}{2\vert y-p\vert}
+\frac{(\vert y-p\vert-T)^3}{3\vert y-p\vert}
\right\},
& y\in B_1\\
\\
\displaystyle
\frac{1}{2}
\left\{
\frac{\eta^3}{6\vert y-p\vert}
-
\frac{\eta(\vert y-p\vert-T)^2}{2\vert y-p\vert}
-\frac{(\vert y-p\vert-T)^3}{3\vert y-p\vert}
\right\},
& y\in B_2
\end{array}
\right.
$$
and
$$
\displaystyle
\partial_tv(y,T)
=
\left\{
\begin{array}{ll}
\displaystyle
\frac{\vert y-p\vert-T}{2\vert y-p\vert}
\left\{\eta-(\vert y-p\vert-T)\right\},
& y\in B_1,\\
\\
\displaystyle
\displaystyle
\frac{\vert y-p\vert-T}{2\vert y-p\vert}
\left\{\eta+(\vert y-p\vert-T)\right\},
& y\in B_2.
\end{array}
\right.
$$

Thus one gets:

(a)  for $y\in B_1$
$$
\displaystyle
v(y,T)
=
\displaystyle
\frac{1}{12}(\eta-2T)(\eta+T)^2\cdot\frac{1}{\vert y-p\vert}
+\frac{1}{2}T(\eta+T)
-\frac{1}{4}(\eta+2T)\vert y-p\vert
+\frac{1}{6}\vert y-p\vert^2
$$
and
$$\displaystyle
\partial_tv(y,T)
=-\frac{1}{2}T(\eta+T)\cdot\frac{1}{\vert y-p\vert}
+\frac{1}{2}(\eta+2T)
-\frac{1}{2}\vert y-p\vert;
$$

(b)  
for $y\in B_2$
$$
\displaystyle
v(y,T)
=
\displaystyle
\frac{1}{12}(\eta+2T)(\eta-T)^2\cdot\frac{1}{\vert y-p\vert}
+\frac{1}{2}T(\eta-T)
-\frac{1}{4}(\eta-2T)\vert y-p\vert
-\frac{1}{6}\vert y-p\vert^2
$$
and
$$\displaystyle
\partial_tv(y,T)
=-\frac{1}{2}T(\eta-T)\cdot\frac{1}{\vert y-p\vert}
+\frac{1}{2}(\eta-2T)
+\frac{1}{2}\vert y-p\vert.
$$

Therefore, we have, for $y\in B_1$,
$$\begin{array}{l}
\displaystyle
\,\,\,\,\,\,
\tau v(y,T)-\partial_tv(y,T)
\\
\\
\displaystyle
=
\left\{\frac{1}{12}\tau(\eta-2T)(\eta+T)^2
+\frac{1}{2}T(\eta+T)\right\}
\cdot\frac{1}{\vert y-p\vert}
+\left\{\frac{1}{2}\tau T(\eta+T)-\frac{1}{2}(\eta+2T)\right\}\\
\\
\displaystyle
\,\,\,
+\left\{-\tau \frac{1}{4}(\eta+2T)+\frac{1}{2}\right\}\vert y-p\vert
+\frac{1}{6}\tau\vert y-p\vert^2,
\end{array}
$$
and, for $y\in B_2$,
$$\begin{array}{l}
\displaystyle
\,\,\,\,\,\,
\tau v(y,T)-\partial_tv(y,T)
\\
\\
\displaystyle
=
\left\{\frac{1}{12}\tau(\eta+2T)(\eta-T)^2
+\frac{1}{2}T(\eta-T)\right\}
\cdot\frac{1}{\vert y-p\vert}
+\left\{\frac{1}{2}\tau T(\eta-T)-\frac{1}{2}(\eta-2T)\right\}\\
\\
\displaystyle
\,\,\,
+\left\{-\tau \frac{1}{4}(\eta-2T)-\frac{1}{2}\right\}\vert y-p\vert
-\frac{1}{6}\tau\vert y-p\vert^2.
\end{array}
$$

Let $x\in B_{T-\eta}(p)\setminus\{p\}$.
Using (3.4) in the case when $R_1=T$, $R_2=T+\eta$, 
we have
$$
\displaystyle
\frac{1}{4\pi}\int_{B_1}\frac{e^{-\tau\vert x-y\vert}}{\vert x-y\vert}
(\tau v(y,T)-\partial_tv(y,T))dy
=\frac{1}{\tau^2}{\cal H}_{+}(\tau;T,\eta)\frac{\sinh\tau\vert x-p\vert}{\vert x-p\vert}.
$$
Using (3.4) in the case when $R_1=T-\eta$, $R_2=T$, 
we have
$$
\displaystyle
\frac{1}{4\pi}\int_{B_2}\frac{e^{-\tau\vert x-y\vert}}{\vert x-y\vert}
(\tau v(y,T)-\partial_tv(y,T))dy
=\frac{1}{\tau^2}{\cal H}_{-}(\tau;T,\eta)\frac{\sinh\tau\vert x-p\vert}{\vert x-p\vert}.
$$
From these we obtain the desired formula.

\noindent
$\Box$

\proclaim{\noindent Proposition 3.3.}
We have
$$\left\{
\begin{array}{l}
\displaystyle
{\cal H}_{+}(\tau;T,\eta)
=f_{\tau}(T)e^{-\tau T}-f_{\tau}(T+\eta)e^{-\tau (T+\eta)},
\\
\\
\displaystyle
{\cal H}_{-}(\tau;T,\eta)
=g_{\tau}(T-\eta)e^{-\tau (T-\eta)}-g_{\tau}(T)e^{-\tau T},
\end{array}
\right.
$$
where
$$\begin{array}{ll}
\displaystyle
f_{\tau}(\xi)
&
\displaystyle
=\frac{\tau}{6}\xi^3+\left\{1-\frac{\tau}{4}(\eta+2T)\right\}\xi^2
+\left\{\frac{1}{2}\tau T(\eta+T)-(\eta+2T)+\frac{2}{\tau}\right\}\xi
\\
\\
\displaystyle
&
\displaystyle
\,\,\,
+\left\{\frac{1}{12}\tau(\eta-2T)(\eta+T)^2
+T(\eta+T)-\frac{\eta+2T}{\tau}+\frac{2}{\tau^2}\right\}
\end{array}
$$
and
$$\begin{array}{ll}
\displaystyle
g_{\tau}(\xi)
&
\displaystyle
=-\frac{\tau}{6}\xi^3-\left\{1+\frac{\tau}{4}(\eta-2T)\right\}\xi^2
+\left\{\frac{1}{2}\tau T(\eta-T)-(\eta-2T)-\frac{2}{\tau}\right\}\xi\\
\\
\displaystyle
&
\displaystyle
\,\,\,
+\left\{\frac{1}{12}\tau(\eta+2T)(\eta-T)^2
+T(\eta-T)-\frac{\eta-2T}{\tau}-\frac{2}{\tau^2}\right\}.
\end{array}
$$

\endproclaim

{\it\noindent Proof.}
First note that we have the relatioship:
$$\left\{
\begin{array}{l}
\displaystyle
H_0(\tau;R_1,R_2)=R_1e^{-\tau R_1}-R_2e^{-\tau R_2}
+\frac{1}{\tau}H_{-1}(\tau;R_1,R_2),
\\
\\
\displaystyle
H_1(\tau;R_1,R_2)
=R_1^2e^{-\tau R_1}-R_2^2e^{-\tau R_2}
+\frac{2}{\tau}H_0(\tau;R_1,R_2),\\
\\
\displaystyle
H_2(\tau;R_1,R_2)
=R_1^3e^{-\tau R_1}-R_2^3e^{-\tau R_2}
+\frac{3}{\tau}H_1(\tau;R_1,R_2).
\end{array}
\right.
\tag {3.7}
$$
Let $R_1=T$ and $R_2=T+\eta$.
Substituting the expression of $H_2(\tau;R_1,R_2)$ in terms of $H_{1}(\tau;R_1,R_2)$ in (3.7) into (3.5),
we have
$$\begin{array}{l}
\displaystyle
\,\,\,\,\,\,
{\cal H}_{+}(\tau;T,\eta)\\
\\
\displaystyle
=
\left\{\frac{1}{12}\tau(\eta-2T)(\eta+T)^2
+\frac{1}{2}T(\eta+T)\right\}H_{-1}(\tau;T,T+\eta)\\
\\
\displaystyle
+\left\{\frac{1}{2}\tau T(\eta+T)-\frac{1}{2}(\eta+2T)\right\}H_{0}(\tau;T,T+\eta)\\
\\
\displaystyle
\,\,\,
+\left\{-\tau \frac{1}{4}(\eta+2T)+\frac{1}{2}\right\}H_{1}(\tau;T,T+\eta)\\
\\
\,\,\,
\displaystyle
+\frac{1}{6}\tau 
\left\{R_1^3e^{-\tau R_1}-R_2^3e^{-\tau R_2}
+\frac{3}{\tau}H_1(\tau;T,T+\eta)\right\}
\\
\\
\displaystyle
=
\left\{\frac{1}{12}\tau(\eta-2T)(\eta+T)^2
+\frac{1}{2}T(\eta+T)\right\}H_{-1}(\tau;T,T+\eta)\\
\\
\displaystyle
\,\,\,
+\left\{\frac{1}{2}\tau T(\eta+T)-\frac{1}{2}(\eta+2T)\right\}H_{0}(\tau;T,T+\eta)\\
\\
\displaystyle
\,\,\,
+\left\{1-\tau \frac{1}{4}(\eta+2T)\right\}H_{1}(\tau;T,T+\eta)\\
\\
\,\,\,
\displaystyle
+\frac{1}{6}\tau(R_1^3e^{-\tau R_1}-R_2^3e^{-\tau R_2}).
\end{array}
$$
Continuing this procedure step by step by using the relationship (3.7) until elliminating all the terms
$H_j(\tau;R_1,R_2)$, $j=1,0$, and finally substituting the explicit form of $H_{-1}(\tau;R_1,R_2)$
into the resulted form,
we obtain 
$$\begin{array}{l}
\,\,\,\,\,\,
{\cal H}_+(\tau;T,\eta)
\\
\\
\displaystyle
=\left\{\frac{1}{12}\tau(\eta-2T)(\eta+T)^2
+T(\eta+T)-\frac{\eta+2T}{\tau}+\frac{2}{\tau^2}\right\}
(e^{-\tau R_1}-e^{-\tau R_2})\\
\\
\displaystyle
\,\,\,
+\left\{\frac{1}{2}\tau T(\eta+T)-(\eta+2T)+\frac{2}{\tau}\right\}
\left(R_1e^{-\tau R_1}-R_2e^{-\tau R_2}\right)
\\
\\
\displaystyle
\,\,\,
+\left\{1-\tau \frac{1}{4}(\eta+2T)\right\}
(R_1^2e^{-\tau R_1}-R_2^2e^{-\tau R_2})
+\frac{1}{6}\tau(R_1^3e^{-\tau R_1}-R_2^3e^{-\tau R_2}).
\end{array}
$$
Making order of this right-hand side, we obtain the dersired expression for ${\cal H}_+(\tau;T,\eta)$.

Next Let $R_1=T-\eta$ and $R_2=T$.  Applying the same procedure based on the relationship (3.7) to the right-hand side on (3.6), we obtain
$$\begin{array}{l}
\displaystyle
\,\,\,\,\,\,
{\cal H}_{-}(\tau;T,\eta)\\
\\
\displaystyle
=\left\{\frac{1}{12}\tau(\eta+2T)(\eta-T)^2
+T(\eta-T)-\frac{\eta-2T}{\tau}-\frac{2}{\tau^2}\right\}
(e^{-\tau R_1}-e^{-\tau R_2})\\
\\
\displaystyle
\,\,\,
+\left\{\frac{1}{2}\tau T(\eta-T)-(\eta-2T)-\frac{2}{\tau}\right\}
\left(
R_1e^{-\tau R_1}-R_2e^{-\tau R_2}
\right)\\
\\
\displaystyle
\,\,\,
+\left\{-1-\tau \frac{1}{4}(\eta-2T)\right\}
(R_1^2e^{-\tau R_1}-R_2^2e^{-\tau R_2})
-\frac{1}{6}\tau(R_1^3e^{-\tau R_1}-R_2^3e^{-\tau R_2}).
\end{array}
$$
This yields the desired expression for ${\cal H}_{-}(\tau;T,\eta)$.

\noindent
$\Box$

From Proposition 3.3 we have
$$\begin{array}{l}
\displaystyle
\,\,\,\,\,\,
{\cal H}_{+}(\tau;T,\eta)+{\cal H}_{-}(\tau;T,\eta)
\\
\\
\displaystyle
=g_{\tau}(T-\eta)e^{-\tau (T-\eta)}+(f_{\tau}(T)-g_{\tau}(T))e^{-\tau T}
-f_{\tau}(T+\eta)e^{-\tau(T+\eta)}.
\end{array}
$$
Moreover, set $\xi=T-\eta$.  Then we have
$$\begin{array}{l}
\displaystyle
\,\,\,\,\,\,
g_{\tau}(T-\eta)
\\
\\
\displaystyle
=-\frac{\tau}{6}\xi^3-\left\{1+\frac{\tau}{4}(-\xi-T)\right\}\xi^2+\left\{-\frac{1}{2}\tau T\xi-(-\xi-T)-\frac{2}{\tau}\right\}\xi\\
\\
\displaystyle
\,\,\,
+\left\{\frac{1}{12}\tau(-\xi+3T)\xi^2
-T\xi-\frac{-\xi-T}{\tau}-\frac{2}{\tau^2}\right\}\\
\\
\displaystyle
=-\frac{\tau}{6}\xi^3
-\xi^2+\frac{\tau}{4}\xi^3+\frac{\tau T}{4}\xi^2
-\frac{\tau T}{2}\xi^2+\xi^2+T\xi-\frac{2\xi}{\tau}\\
\\
\displaystyle
\,\,\,
-\frac{\tau}{12}\xi^3+\frac{\tau T}{4}\xi^2
-T\xi+\frac{\xi}{\tau}+\frac{T}{\tau}-\frac{2}{\tau^2}\\
\\
\displaystyle
=\frac{T}{\tau}-\frac{\xi}{\tau}-\frac{2}{\tau^2}\\
\\
\displaystyle
=\frac{\eta}{\tau}-\frac{2}{\tau^2}\\
\\
\displaystyle
=\frac{1}{\tau}\left(\eta-\frac{2}{\tau}\right).
\end{array}
\tag {3.8}
$$

\newpage

\noindent
This yields
$$\begin{array}{l}
\displaystyle
\,\,\,\,\,\,
e^{\tau(T-\eta)}(
{\cal H}_{+}(\tau;T,\eta)+{\cal H}_{-}(\tau;T,\eta))
\\
\\
\displaystyle
=g_{\tau}(T-\eta)+O(\tau e^{-\tau\eta})
\\
\\
\displaystyle
=\frac{1}{\tau}(\eta+O(\tau^{-1})).
\end{array}
$$
This completes the proof of Lemma 2.2.

{\bf\noindent Remark 3.2.}  Similar to the derivation of (3.8), one gets
$$
\left\{
\begin{array}{l}
\displaystyle
f_{\tau}(T)=\frac{\tau}{12}\eta^3-\frac{\eta}{\tau}+\frac{2}{\tau^2},\\
\\
\displaystyle
f_{\tau}(T+\eta)=\frac{1}{\tau}\left(\eta+\frac{2}{\tau}\right),\\
\\
\displaystyle
g_{\tau}(T)=\frac{\tau}{12}\eta^3-\frac{\eta}{\tau}-\frac{2}{\tau^2},
\end{array}
\right.
$$
and thus
$$\begin{array}{l}
\displaystyle
\,\,\,\,\,\,
{\cal H}_+(\tau,T,\eta)
+{\cal H}_{-}(\tau;T,\eta)
\\
\\
\displaystyle
=\frac{4}{\tau^2}e^{-\tau T}
-\frac{1}{\tau}\left(\eta+\frac{2}{\tau}\right)e^{-\tau(T+\eta)}
+\frac{1}{\tau}\left(\eta-\frac{2}{\tau}\right)e^{-\tau(T-\eta)}.
\end{array}
$$
However, we do not need this explicit formula for the present purpose.

$$\quad$$

\centerline{{\bf Acknowledgments}}

The author would like to thank anonymous referees for giving 
valuable comments on the improvement of the presentation of the results.
The author was partially supported by Grant-in-Aid for
Scientific Research (C)(No. 17K05331) and (B) (No. 18H01126) of Japan  Society for
the Promotion of Science.

$$\quad$$

\section{Appendix}

In this appendix we give an explicit computation result for the potential
$$\displaystyle
v_j(x)=\int_B\frac{e^{-\tau\vert x-y\vert}}{\vert x-y\vert}\vert y\vert^jdy,\,\,x\in B,
$$
where $B=\{y\in\Bbb R^3\,\vert\,\vert y\vert<\eta\}$, with $\eta>0$ and $j=-1,0,1,2$.

\proclaim{\noindent Proposition A.}
For all $x\in B\setminus\{0\}$ we have
$$\displaystyle
\left\{
\begin{array}{l}
\displaystyle
v_{-1}(x)
=\frac{4\pi}{\tau^2}
\left(\frac{1-e^{-\tau\vert x\vert}}{\vert x\vert}-e^{-\tau\eta}\frac{\sinh\tau\vert x\vert}{\vert x\vert}\right),
\\
\\
\displaystyle
v_0(x)=\frac{4\pi}{\tau^2}
\left\{
1-\left(\eta+\frac{1}{\tau}\right)e^{-\tau\eta}\frac{\sinh\tau\vert x\vert}{\vert x\vert}\right\},
\\
\\
\displaystyle
v_{1}(x)=\frac{4\pi}{\tau^2}\left\{\vert x\vert+\frac{2}{\tau^2}\frac{1-e^{-\tau\vert x\vert}}{\vert x\vert}
-
e^{-\tau\eta}\left(\eta^2+\frac{2}{\tau}
\eta+\frac{2}{\tau^2}\right)\frac{\sinh\tau\vert x\vert}{\vert x\vert}
\right\},
\\
\\
\displaystyle
v_2(x)=
\frac{4\pi}{\tau^2}
\left\{
\vert x\vert^2+\frac{6}{\tau^2}-e^{-\tau\eta}
\left(\eta^3
+\frac{3\eta^2}{\tau}+\frac{6\eta}{\tau^2}+\frac{6}{\tau^3}\right)\frac{\sinh\tau\vert x\vert}{\vert x\vert}
\right\}.
\end{array}
\right.
$$

\endproclaim

{\it\noindent Proof.}
The change of variables $y=r\omega\,(0<r<\eta,\,\omega\in S^2)$ and a rotation give us
$$\begin{array}{ll}
\displaystyle
v_j(x)
&
\displaystyle
=\int_0^{\eta}r^{2+j} dr\int_{S^2}\frac{e^{-\tau\vert x-r\omega\vert}}{\vert x-r\omega\vert}d\omega
\\
\\
\displaystyle
&
\displaystyle
=\int_0^{\eta}r^{2+j} dr\int_{S^2}
\frac{\displaystyle e^{-\tau\vert \vert x\vert\mbox{\boldmath $e$}_3-r\omega\vert}}
{\displaystyle\vert\vert x\vert\mbox{\boldmath $e$}_3-r\omega\vert}d\omega\\
\\
\displaystyle
&
\displaystyle
=\int_0^{\eta}r^{2+j}dr\int_0^{2\pi}d\theta
\int_0^{\pi}\sin\varphi d\varphi
\frac{\displaystyle e^{-\tau\sqrt{\vert x\vert^2-2r\vert x\vert\cos\varphi+r^2}}}
{\displaystyle\sqrt{\vert x\vert^2-2r\vert x\vert\cos\varphi+r^2}}\\
\\
\displaystyle
&
\displaystyle
=2\pi\int_0^{\eta}Q(\vert x\vert,r)r^{2+j} dr,
\end{array}
$$
where $\mbox{\boldmath $e$}_3=(0,0,1)$ and 
$$\displaystyle
Q(\xi,r)
=\int_0^{\pi}
\frac{\displaystyle e^{-\tau\sqrt{\xi^2-2r\xi\cos\varphi+r^2}}}
{\displaystyle\sqrt{\xi^2-2r\xi\cos\varphi+r^2}}\sin\varphi d\varphi,\,0\le\xi<\eta,\,0<r<\eta.
$$
Fix $\xi\in]0,\,\eta[$ and $r\in]0,\,\eta[$.
The change of variable
$$\displaystyle
s=\sqrt{\xi^2-2r\xi\cos\varphi+r^2},\,\varphi\in]0,\,\pi[,
$$
gives
$$\displaystyle
s^2=\xi^2-2r\xi\cos\varphi+r^2
$$
and
$$\displaystyle
sds=r\xi\sin\varphi d\varphi.
$$
Hence, we have
$$\begin{array}{ll}
\displaystyle
Q(\xi,r)
&
\displaystyle
=\frac{1}{r\xi}\int_{\vert\xi-r\vert}^{\xi+r}e^{-\tau s}ds
\\
\\
\displaystyle
&
\displaystyle
=-\frac{1}{r\xi\tau}\left(e^{-\tau(\xi+r)}-e^{-\tau\vert\xi-r\vert}\right).
\end{array}
$$
Therfore, we obtain
$$
\begin{array}{ll}
\displaystyle
v_j(x) &
\displaystyle
=2\pi\int_0^{\eta}Q(\vert x\vert,r)r^{2+j} dr
\\
\\
\displaystyle
&
\displaystyle
=\frac{2\pi}{\xi\tau}\int_0^{\eta}\left(e^{-\tau\vert\xi-r\vert}-e^{-\tau(\xi+r)}\right)r^{1+j}dr\vert_{\xi=\vert x\vert}.
\end{array}
\tag {A.1}
$$
Thus, everything is reduced to computing the integral
$$\begin{array}{ll}
\displaystyle
K_j=
\int_0^{\eta}\left(e^{-\tau\vert\xi-r\vert}-e^{-\tau(\xi+r)}\right)r^{1+j}\,dr, & j=-1,0,1,2.
\end{array}
$$

A direct computation yields
$$
\displaystyle
\left\{
\begin{array}{l}
\displaystyle
\int_0^{\eta}e^{-\tau\vert\xi-r\vert}dr
=\frac{2}{\tau}-\frac{e^{-\tau\xi}}{\tau}-\frac{e^{-\tau(\eta-\xi)}}{\tau},
\\
\\
\displaystyle
\int_0^{\eta}e^{-\tau\vert\xi-r\vert}rdr
=\frac{2\xi}{\tau}
+\frac{e^{-\tau\xi}}{\tau^2}
-\frac{e^{\tau(\xi-\eta)}}{\tau}\left(\eta+\frac{1}{\tau}\right),\\
\\
\displaystyle
\int_0^{\eta}e^{-\tau\vert\xi-r\vert}r^2dr
=\frac{1}{\tau^3}
\left\{(2\tau^2\xi^2+4)-2e^{-\tau\xi}
-(\tau^2\eta^2+2\tau\eta+2)e^{-\tau(\eta-\xi)}
\right\},
\\
\\
\displaystyle
\int_0^{\eta}e^{-\tau\vert\xi-r\vert}r^3dr
=\frac{2\xi^3}{\tau}-\frac{1}{\tau} e^{-\tau(\eta-\xi)}\eta^3
+\frac{6}{\tau^4}e^{-\tau\xi}
-\frac{3}{\tau^4}
\left\{
(\tau^2\eta^2+2\tau\eta+2)e^{-\tau(\eta-\xi)}
-4\tau\xi
\right\}.
\end{array}
\right.
$$
Also we have
$$
\displaystyle
\left\{
\begin{array}{l}
\displaystyle
\int_0^{\eta}e^{-\tau(\xi+r)}dr
=\frac{e^{-\tau\xi}}{\tau}-\frac{e^{-\tau(\xi+\eta)}}{\tau},
\\
\\
\displaystyle
\int_0^{\eta}e^{-\tau(\xi+r)}rdr
=\frac{e^{-\tau\xi}}{\tau^2}
-\frac{e^{-\tau(\xi+\eta)}}{\tau}
\left(\eta+\frac{1}{\tau}\right),
\\
\\
\displaystyle
\int_0^{\eta}e^{-\tau(\xi+r)}r^2dr
=\frac{1}{\tau^3}
e^{-\tau\xi}\left\{-e^{-\tau\eta}(\tau^2\eta^2+2\tau\eta+2)+2\right\},
\\
\\
\displaystyle
\int_0^{\eta}e^{-\tau(\xi+r)}r^3dr
=-\frac{1}{\tau}\eta^3e^{-\tau(\xi+\eta)}
+\frac{3}{\tau^4}e^{-\tau\xi}
\left\{-e^{-\tau\eta}(\tau^2\eta^2+2\tau\eta+2)+2\right\}.
\end{array}
\right.
$$
From these, we obtain
$$\displaystyle
\left\{
\begin{array}{l}
\displaystyle
K_{-1}
=\frac{2}{\tau}(1-e^{-\tau\xi}-e^{-\tau\eta}\sinh\tau\xi),\\
\\
\displaystyle
K_0
=\frac{2}{\tau}
\left\{\xi-\left(\eta+\frac{1}{\tau}\right)e^{-\tau\eta}\sinh\tau\xi\right\},
\\
\\
\displaystyle
K_1
=\frac{2}{\tau^3}
\left\{(\tau^2\xi^2+2)-2e^{-\tau\xi}-(\tau^2\eta^2+2\tau\eta+2)e^{-\tau\eta}\sinh\tau\xi
\right\},
\\
\\
\displaystyle
K_2
=
\frac{2\xi^3}{\tau}+\frac{12}{\tau^3}\xi
-\frac{2}{\tau}\left(\eta^3
+\frac{3\eta^2}{\tau}+\frac{6\eta}{\tau^2}+\frac{6}{\tau^3}\right)e^{-\tau\eta}\sinh\tau\xi.
\end{array}
\right.
$$
Substituting these into (A.1), we obtain the desired formulae.

\noindent
$\Box$

\vskip1cm
\noindent
e-mail address

ikehata@hiroshima-u.ac.jp


\begin{thebibliography}{99}






\bibitem{BEL}  Belishev, M. I.,
             On an approach to multidimensional inverse problems for the wave equation,
             Dokl. Akad. Nauk SSSR, {\bf 297}(1987), 524-527.




\bibitem{B2} Belishev, M. I.,
             How to see waves under the Earth surface (the BC-method for geophysicists),
             Ill-Posed and Inverse Problems, pp. 67-84,
             Kabanikhin, S. I. and Romanov, V. G. (Eds), VSP, Utrecht, 2002.






\bibitem{BKLS} Bingham, K., Kurylev, Y., Lassas, M. and Siltanen, S.,
              Iterative time-reversal control for inverse problems, 
              Inverse Problems and Imaging, {\bf 2}(2008), 63-81.







\bibitem{BP}  Burkard, C. and Potthast, R.,
              A time-domain probe method for three-dimensional rough surface reconstructions,
              Inverse Problems and Imaging, {\bf 3}(2009), pp. 259-274.






\bibitem{DKL} Dahl, M. F., Kirpichnikova, A. and Lassas, M.,
              Focusing waves in unknown media by modified time reversal iteration,
              SIAM J. Control. Optim., {\bf 48}(2009), 839-858.










\bibitem{dBK} de Buhan, M. and Kray, M., 
               A new approach to solve the inverse scattering problem for waves: 
               combining the TRAC and the Adaptive Inversion methods,
               Inverse Problems, {\bf 29}(2013), 085009.
               
               
               


























\bibitem{Du}  Duff,  G.F.D., 
              Hyperbolic Differential Equations and Waves. 
              In: Garnir H.G. (eds) Boundary Value Problems for Linear Evolution Partial Differential Equations, pp.27-155,
              NATO Advanced Study Institutes Series (Series C-Mathematical and Physical Sciences), vol 29. Springer, Dordrecht, 1977.





\bibitem{FWCM} Fink, M.,
               Time reversal of ultrasonic fields-Part I: Basic principles,
               IEEE Trans. Ultrason., Ferroelec., Freq. Contr., {\bf 39}(1992), No.5, pp.555-566.













\bibitem{I}  Ikawa, M., Mixed problems for hyperbolic equations of second order, 
             J. Math. Soc. Japan, {\bf 20}(1968), 580-608.








\bibitem{E00} Ikehata, M.,
              \newblock Reconstruction of the support function for inclusion from boundary measurements,
              J. Inverse Ill-Posed Problems, {\bf 8}(2000), 367-378.















\bibitem{IW00} Ikehata, M.,
               The enclosure method for inverse obstacle scattering problems with dynamical data over a
               finite time interval, Inverse Problems, {\bf 26}(2010) 055010(20pp).








\bibitem{IEO2} Ikehata, M.,
              The enclosure method for inverse obstacle scattering problems with dynamical data over a finite time
              interval: II. Obstacles with a dissipative boundary or finite refractive index and back-scattering data,
              Inverse Problems, {\bf 28}(2012) 045010 (29pp).






\bibitem{IEO3} Ikehata, M.,
               The enclosure method for inverse obstacle scattering problems with dynamical data over a finite time interval:
               III. Sound-soft obstacle and bistatic data, Inverse Problems, {\bf 29}(2013) 085013 (35pp).




























\bibitem{EIV} Ikehata, M., 
              The enclosure method for inverse obstacle scattering over a finite time interval: IV.
              Extraction from a single point on the graph of the response operator,
              J. Inverse Ill-Posed Probl., {\bf 25}(2017), 747-761.
              


        






        
 

\bibitem{II}  Ikehata, M. and Itou, H., On reconstruction of a cavity in a linearized viscoelastic body from infinitely many transient
              boundary data, Inverse Problems, {\bf 28}(2012) 125003 (19pp).








\bibitem{IK1}  Ikehata, M. and Kawashita, M., The enclosure method for the heat equation,
               Inverse Problems, {\bf 25}(2009) 075005(10pp).
               
               
               



\bibitem{IK2}  Ikehata, M. and Kawashita, M., On the reconstruction of inclusions in a heat conductive body
               from dynamical boundary data over a finite time interval, Inverse Problems,
               {\bf 26}(2010) 095004(15pp).
    












\bibitem{Is} Isakov, V., 
               Inverse obstacle problems, TOPICAL REVIEW, Inverse Problems, {\bf 25}(2009) 123002(18pp).
    
          





\bibitem{LCW}  Lines, C. D. and Chandler-Wilde, S. N.,
               A time domain point source method for inverse scattering by rough surfaces,
               Computing, {\bf 75}(2005), No. 2, pp. 157-180.






\bibitem{MS}  Monk, P. and Selgas, V.,
              An inverse acoustic waveguide problem in the time domain,
              Inverse Problems, {\bf 32}(2016), 05501(26pp).








\bibitem{O} Oksanen, L., 
             Solving an inverse obstacle problem for the wave equation by using the boundary control method,
             Inverse Problems, {\bf 29}(2013) 035004.











\bibitem{Y}   Yosida, K., Functional Analysis, Third Edition, Springer, New York, 1971.




\end{thebibliography}
\end{document}